\documentclass[reqno,12pt]{amsart}
\usepackage{amsfonts} %如果在方括号里面加上draft，可以看到文中是否有过长的句子。
\usepackage{times}
\usepackage{amssymb, amsmath}
\allowdisplaybreaks

\begin{document}

\title[ Inviscid Limit for Vortex Patches]{Inviscid Limit for Vortex Patches\\[3mm] in A Bounded Domain}

\author[Jiu,Wang]{}

\maketitle

 \centerline{\scshape Quansen Jiu\footnote{The research is partially
supported by National Natural Sciences Foundation of China (No.
10871133) and Project of Beijing Education Committee.}}
\medskip
{\footnotesize
  \centerline{School of Mathematical Sciences, Capital Normal University}
  \centerline{Beijing  100048, P. R. China}
   \centerline{  \it Email: qsjiumath@gmail.com}
   }
\medskip
\centerline{\scshape Yun Wang}
\medskip
{\footnotesize
  \centerline{The Institute of Mathematical Sciences,
  The Chinese University of Hong Kong}
  \centerline{Shatin, N.T., Hong Kong}
   \centerline{  \it Email: ywang@math.cuhk.edu.hk }
   }

\begin{center}
\begin{minipage}{13cm}

Abstract: In this paper, we consider the inviscid limit  of
 the incompressible Navier-Stokes equations  in a
smooth, bounded and simply connected domain $\Omega \subset
\mathbb{R}^d, d=2,3$.  We prove that for a vortex patch initial data
the weak Leray solutions of the incompressible Navier-Stokes
equations with Navier boundary conditions will converge (locally in
time for $d=3$ and globally in time for $d=2$) to a vortex patch
solution of the incompressible Euler equation as the viscosity
vanishes. In view of the results obtained in \cite{ad} and \cite{mo}
which dealt with the case of the whole space, we derive an almost
optimal convergence rate $(\nu t)^{\frac34-\varepsilon}$ for any
small $\varepsilon>0$ in $L^2$.
\end{minipage}
\end{center}

\vspace{3mm}\textbf{Keywords}:~~inviscid limit, Navier boundary
condition, vortex patches

\vspace{2mm}\textbf{AMS Subject Classifications}: 35Q~~76D

\vspace{5mm} \section{Introduction}

 \setcounter{equation}{0}

\vspace{2mm} The incompressible Navier-Stokes equations read as
\begin{equation}\label{1.1}
\left\{\begin{array}{l} \displaystyle\frac{\partial u}{\partial t}
-\nu\Delta u + (u\cdot\nabla)u + \nabla p  = 0,\\
[2mm] {\rm div}~u = 0,
\end{array}\right.
\end{equation}
where $u=(u_1,\cdots, u_d) (d=2$  or $3)$ is the velocity fields,
$p$ is the pressure function and $\nu$ is the kinetic viscosity.

Formally, when $\nu=0$, \eqref{1.1} becomes the following
incompressible Euler equations:
\begin{equation}\label{1.1+}
\left\{\begin{array}{l} \displaystyle\frac{\partial u}{\partial t}
 + (u\cdot\nabla)u + \nabla p  = 0,\\
[2mm] {\rm div}~u = 0.
\end{array}\right.
\end{equation}

The inviscid limit for the incompressible Navier-Stokes equations in
the whole space  has been well understood (see
\cite{Swan,che,mo,P1,ad,cw1,cw2} and references therein) for both
smooth and non-smooth initial data. However,  in the case of a
bounded domain, the inviscid limit for the Navier-Stokes equations
 with Dirichlet boundary conditions is still a
completely open problem. This is mainly due to the difference
between the Dirichlet boundary conditions of the incompressible
Navier-Stokes equations \eqref{1.1} and the tangential boundary
conditions
 of the incompressible Euler equations \eqref{1.1+} and  a boundary
 layer will appear near the boundary of the domain.

This paper is concerned with  the inviscid limit problem of the
incompressible Navier-Stokes equations \eqref{1.1} with the
following Navier boundary conditions:
\begin{equation}\label{1.2}
u\cdot \vec{n}=0,\ \ \ \ [D(u)\vec{n} + \alpha u]_{tan} =0, \ \ \
on\
\partial \Omega \times (0, +\infty),
\end{equation}
where $\Omega\subset\mathbb{R}^d (d=2$ or $3)$ is a smooth bounded
domain, $\vec{n}$ is the unit exterior normal to the boundary
$\partial\Omega$, $D(u) = \frac{1}{2} [\nabla u + (\nabla u)^T] $ is
the rate of strain tensor and $[D(u)\vec{n} + \alpha u]_{tan}$ is
the tangential component of the vector $D(u)\vec{n} + \alpha u$.
Here $\alpha=\alpha(x,t)$ is a known function representing the
friction coefficient of the material.

 The Navier boundary conditions,  introduced by Navier in \cite{na}, say that the
tangential component of the viscous stress at the boundary is
proportional to the tangential velocity. They were rigorously
justified as a homogenization of the no-slip condition on a rough
boundary in \cite{JM} and  widely used when studying the inviscid
limit of the incompressible flows in a bounded domain (see
\cite{cmr, ip, JN, ke, lnp,X1}) in recent years.

Of particular interest of this paper is the inviscid limit for
vortex patches in a bounded domain.  It is known that when the
initial data are vortex patch ones (see Definition 2.1 for details),
there exists a unique solution to the incompressible Euler equations
which preserves the vortex-patch structures globally (in time) in
the whole plane (\cite{che1,bc}) and locally in three-dimensional
whole space ( \cite{GS}).  In a smooth, bounded and simply connected
domain $\Omega\subset \mathbb{R}^d (d=2,3)$, the vortex patch
solutions of the incompressible Euler equations were derived in
\cite{du,de}. In this paper, we will show that the  weak Leray
solutions of the incompressible Navier-Stokes equations with Navier
boundary conditions  will tend to a vortex patch solution of the
incompressible
 Euler equations as the viscosity vanishes if the initial data are vortex patch ones.
  Moreover, we obtain that the convergence rate in $L^2$ is $(\nu t)^{\frac34-\varepsilon}$
  for any small $\varepsilon>0$. In the case of the whole plane,
 Constantin and Wu
studied the inviscid limit for the 2D vortex patches  in
\cite{cw1,cw2} and obtained the convergence rate in $L^2$ is
$\sqrt{\nu t}$. Abidi and Danchin  improved the convergence rate  to
be $(\nu t)^{\frac34}$ in $L^2$  which is  optimal since the
circular vortex patches provide a lower bound (see \cite{ad}).
Later, Masmoudi extended the results  to the case of
three-dimensional whole space  in \cite{mo}. Recently, Sueur
\cite{S} dealt  with the vorticity internal transition layers
 for the Navier-Stokes equations and described how the smoothing effect is (micro-)localized in the case
where vortex patches are prescribed as initial data, using the
method of asymptotic expansion. In the case of the two-dimensional
bounded domain, the inviscid limit for the incompressible
Navier-Stokes equations with Navier- boundary conditions was
discussed in \cite{ke} and the obtained convergence rate in $L^2$ is
$\sqrt{\nu t}$ for initial vorticity in $L^\infty$.

Our results here applies to both 2D and 3D vortex patches in a
bounded domain and the convergence rate obtained in this paper is
almost optimal in view of the results in \cite{ad} and \cite{mo}.
Since we consider the case  of the bounded domain, estimates in
Besov space in \cite{ad,mo} can not be used directly and we will use
the interpolation space theory to deduce that the vorticity belongs
to $L^\infty([0,T^*);H^s(\Omega)$ for some $T^*>0$ and $s>0$. More
subtle estimates will be given in this paper. Meanwhile, whether the
convergence rate can be improved to $(\nu t)^{\frac34}$ is still
open.

 The paper is organized as follows. In
Section 2, we will give some preliminaries and the main results.
Section 3 is devoted to the proof of the main result.

%%%%%%%%%%%%%%%%%%% Main result %%%%%%%%%%%%%%%%%%%%%%%%%%%%%%%%%%%%%%%%%%%%%%%%%%%%%
 \vspace{5mm}

\section{Preliminaries and Main Results}
\setcounter{equation}{0} \vspace{2mm}

Let  $\Omega \subset \mathbb{R}^d (d=2,3)$ be a smooth, bounded and
simply connected domain. The initial-boundary problem to the
incompressible Euler equations is written as
\begin{equation}\label{1.3}
\left\{\begin{array}{l} \displaystyle\frac{\partial u}{\partial t}
 + (u\cdot\nabla)u + \nabla p  = 0, \ \ \ \ \ (x,t)\in \Omega\times (0,+\infty)\\
[2mm] {\rm div}~u = 0, \ \ \ (x,t)\in \Omega \times [0, +\infty),\\[1mm]
u\cdot \vec{n} =0, \ \ \ \ \ (x,t)\in \partial \Omega\times
[0,+\infty),\\[2mm]
u(x,0)=u_0(x),\ \ \ \ \ x\in \Omega.
\end{array}\right.
\end{equation}

Denote by $(u^\nu, p^\nu)$ the solutions of the incompressible
Navier-Stokes equations with corresponding kinetic viscosity $\nu$.
The initial-boundary problem to the incompressible Navier-Stokes
equations with Navier boundary conditions is written as
\begin{equation}\label{1.3+}
\left\{\begin{array}{l} \displaystyle\frac{\partial u^\nu}{\partial
t}-\nu\Delta u^\nu
 + (u^\nu\cdot\nabla)u^\nu + \nabla p^\nu  = 0, \ \ \ \ \ (x,t)\in \Omega\times (0,+\infty)\\
[2mm] {\rm div}~u^\nu = 0, \ \ \ (x,t)\in \Omega \times [0, +\infty),\\[1mm]
u^\nu\cdot \vec{n}=0,\ \ \ \ [D(u^\nu)\vec{n} + \alpha u^\nu]_{tan}
=0, \ \ \ (x,t)\in
\partial \Omega \times
[0,+\infty),\\[2mm]
u(x,0)=u^\nu_0(x),\ \ \ \ \ x\in \Omega.
\end{array}\right.
\end{equation}

 Let $\omega_0={\rm curl}~u_0$ be the initial vorticity of $u_0$.
In this paper, for any vector-valued function $\varphi$,
$D(\varphi)$ denotes the symmetric part of $\nabla \varphi$, i.e.,
$$D(\varphi) = \frac{\nabla \varphi + (\nabla \varphi)^T}{2}.$$

 Denote by $C^r, C^{1+r}$
 ($0<r<1 $) the usual H\"older space. In particular, $C_c^r(\mathbb{R}^d)$ consists of functions
 in $C^r(\mathbb{R}^d)$ with compact support.
Let $L^p(\Omega), W^{s,p}(\Omega)$ be the usual Sobolev spaces
defined in $\Omega$, where $1\le p\le \infty$ and $s$ is permitted
to be a real number. If $p=2$,  $W^{s,2}(\Omega)$ is denoted by
$H^s(\Omega)$. $H_0^s(\Omega)$ is the closure of
$C_0^\infty(\Omega)$ in $H^s(\Omega)$. Define
$$
C^\infty_{0,\sigma}(\Omega) =\{f| f\in C_0^\infty(\Omega),\  {\rm
div}~f=0\},
$$
$$C^\infty_{\sigma}(\Omega) =\{f| f\in C^\infty(\overline{\Omega}),\  {\rm
div}~f=0\}.$$  $L^2_\sigma(\Omega)$  is  the closure of
$C^\infty_{0,\sigma}(\Omega)$ in $L^2(\Omega)$, and
$H^1_\sigma(\Omega)$ is  the closure of $C^\infty_{\sigma}(\Omega)$
in $H^1(\Omega)$.

\vspace{2mm}

We first recall the definition of a vortex patch in a bounded domain
(see \cite{du, mo}). \vspace{1mm}

\textbf{Definition 2.1}~~Let $0<r<1$. The vorticity $\omega={\rm
curl}~u$ of a vector field $u$ is called a $C^r$  vortex patch of
support $P$ if  the following decomposition holds:
$$\omega =( \omega_{i}\chi_{P} + \omega_{e}\chi_{\Omega\setminus
\overline{P}})|_\Omega,$$ where $P$ is an open set of class
$C^{1+r}$, $\omega_{i}, \omega_{e} \in C_c^r(\mathbb{R}^d) (d=2,3)$
and $\chi_{P}$, $\chi_{\Omega\setminus \overline{P}}$ are the
characteristic functions of $P$ and $\Omega \setminus \overline{P}$
respectively.

Notice that when $d=3$, ${\rm curl}~u$ is of divergence free, we
need $\omega_i\cdot\vec{n}=\omega_e\cdot\vec{n}$ on $\partial P$.

 If the initial data of the incompressible Euler equations is a  $C^r$ vortex
patch, the global existence of 2-d vortex patch solutions and the
local existence of 3-d vortex patch solutions have been proved (see
\cite{du,de}). More precisely, one has\vspace{2mm}

\textbf{Theorem 2.1}~~ Let $u_0$ be a divergence free vector field
in $\mathbb{R}^d (d=2,3)$, tangent to $\partial \Omega$, whose
vorticity $\omega_0$ is a $C^r$ vortex patch of support $P$, the
boundary of $\partial P$ is a $(d-1)$-dimensional compact
submanifold of $\mathbb{R}^d$. If $\overline{P}\subset \Omega$, then
there exists a $T^*>0$ such that the Euler equations \eqref{1.3}
have
 a (unique) solution $u\in L^\infty([0,T^*); Lip(\Omega)).$ Moreover, $\omega(t) ={\rm curl}~u(t)$
 remains a vortex patch, whose
support $\Psi(t,P)$ is of class $C^{1+r}$ for any $t\in [0,T^*), $
$\Psi$ denoting the flow of $u$. In addition, $T^*>0$ can be
arbitrarily large if $d=2$.

\vspace{2mm}\textbf{Remark 2.1} Under assumptions of Theorem 2.1, if
$P$ is tangent to $\partial\Omega$, a little regularity may be lost.
However, local existence of 3-D vortex patch of $C^s$$(0<s<r)$ and
global existence of 2-D vortex patch of $C^s$$(0<s<r)$ is proved in
\cite{du}.

\vspace{2mm}

 Now we give the definition of a Leray weak
solution of the incompressible Navier-Stokes equations with Navier
boundary conditions.

 \vspace{2mm} \textbf{Definition 2.2}~~We call a
 vector field $u^\nu(t,x): [0,+\infty) \times \Omega
\rightarrow \Omega $, denoted by $u(t,x)$, a weak Leray solution of
\eqref{1.3+} if $u$ verifies
\begin{enumerate}
\item $ u \in C_w([0,\infty); L_{\sigma}^2(\Omega))\bigcap L_{loc}([0,\infty);
H_\sigma^1(\Omega));$

\item $u$ verifies the system of equations \eqref{1.3+}
under the following weak form: for every $\varphi \in C_0^\infty([0,
\infty); C^\infty_\sigma(\Omega))$ with $\varphi\cdot \vec{n}=0$ on
$\partial \Omega$,
$$
\begin{array}{ll} &\displaystyle
2\alpha \nu \int_0^\infty \int_{\partial\Omega} u\cdot \varphi +
2\nu \int_0^\infty \int_\Omega D(u)D(\varphi) + \int_0^\infty
\int_\Omega (u\cdot\nabla)u\cdot \varphi\\[3mm]
& =\displaystyle \int_0^\infty \int_\Omega u\cdot \partial_t \varphi
+\int_\Omega u(0) \cdot \varphi(0) \end{array}$$

\item $u$ verifies the energy inequality, for all $t\geq 0$,
$$\|u(t)\|_{L^2(\Omega)}^2 + 4\alpha \nu \int_0^t \int_{\partial
\Omega} |u|^2 + 4\nu \int_0^t \int_\Omega|D(u)|^2 \leq
\|u(0)\|_{L^2(\Omega)}^2.$$
\end{enumerate}

We remark that the global existence of the Leray weak solution in
the case of Direchlet boundary conditions is well known for any
$u_0\in L^2_\sigma(\Omega)$. The extensions of this result to the
case of Navier boundary conditions is straightforward by the
Galerkin method.

\vspace{2mm}

 The main result of the paper is stated as

\vspace{2mm}\textbf{Theorem 2.2}~~Suppose that the assumptions of
Theorem 2.1 hold and  $u\in L^\infty([0, T^*);
Lip(\overline{\Omega}))$ is the vortex patch solution of the
incompressible Euler equations with initial data $u_0$. Suppose that
$u^\nu$ are Leray weak solutions of the incompressible Navier-Stokes
equations with Navier boundary conditions \eqref{1.3+}. The
corresponding initial data $u^\nu(0)$ is uniformly bounded in
$L^2(\Omega)$, and $\alpha \in L^\infty(\partial \Omega).$ Then for
all $0< T < T^*$ and any small $\epsilon>0$, one has
$$\|(u^\nu - u)(t)\|_{L^2(\Omega)} \leq C ( (\nu t)^{\frac{1+\beta - \epsilon}{2} }+
\|u^{\nu}(0) - u_0\|_{L^2(\Omega)}),$$ where $\beta =
\min(\frac{1}{2}, r)$ and $C$ is a constant depending only on
$\epsilon, u, T,$ $\|\alpha\|_{L^\infty(\partial \Omega)}$ and
$M\equiv \sup_\nu \|u^\nu(0)\|_{L^2(\Omega)}$.

%%%%%%%%%%%%%%%%%%%%%%%%%%%%%%%%%%%%%Proof%%%%%%%%%%%%%%%%%%%%%%%%%%%%%%%
 \vspace{5mm}

\section{Proof of Main Result}

Since we are concerned with the case of the bounded domain, the
estimates in Besov space as in \cite{ad} and \cite{mo} can not be
used directly. However, we have

\vspace{2mm} \textbf{Lemma 3.1}~~Suppose that $\omega = {\rm
curl}~u$ is the vortex solution to the incompressible Euler system,
derived in Theorem 2.1. Then, for any $s <\beta = \min
(r,\frac{1}{2}),$ one has $\omega \in L^\infty ([0, T^*);
H^s(\Omega)).$

\textbf{Proof.}~~It is proved in \cite{du} that the vortex patch
solution has the following structures:
$$\omega(x,t) = \omega_i(x,t)\chi_{P(t)}(x) + \omega_e(x,t)\chi_{\Omega
\setminus \overline{P(t)}}(x),   t\in [0,T^*),$$ where
$$\omega_i,
\omega_e \in L^\infty([0, T^*); C^{\tilde{r}}(\mathbb{R}^d),\ \
P(t)\in L^{\infty}([0,T^*);C^{1+\tilde{r}}(\mathbb{R}^d))$$ for any
$\tilde{r} < r$, which means that for any $ t\in [0,T^*), $ $P(t)$
is a $C^{1+\tilde{r}}$ domain, and the $C^{1+\tilde{r}}$-norm of the
boundary $\partial P(t)$ is locally bounded. Hence
$\mathcal{H}^{d-1}(\partial P(t))$, the $(d-1)$-dimensional
Hausdorff measure of $\partial P(t)$ is locally bounded which
induces that
$$\chi_{P(t)}(x),\ \chi_{\Omega\setminus \overline{P(t)}}(x) \in
L^{\infty} ([0,T^*); L^\infty (\mathbb{R}^d) \cap
BV(\mathbb{R}^d)).$$ Following Lemma 4.2 and Lemma 4.3 in \cite{mo},
after extending $\omega$ to the whole space by zero extension, we
derive
$$\omega(x,t) \in L^\infty ([0, T^*);
\dot{B}^{\tilde{s}}_{2,\infty}(\mathbb{R}^d)),$$ where $\tilde{s} =
\min (\tilde{r}, \frac{1}{2})$ and $\dot{B}_{2,\infty}^{\tilde{s}} $
is the classical homogeneous Besov space (see \cite{tr} for
definition).

Using the fact that $\omega(x,t) \in L^\infty([0, T^*);
L^2(\mathbb{R}^d))$, one has
\begin{equation}\label{9-11-1}
\omega(x,t) \in L^\infty ([0,T^*);
\dot{B}_{2,\infty}^{\tilde{s}}(\mathbb{R}^d) \cap
L^2(\mathbb{R}^d)).
\end{equation}
Moreover, for any $s<\beta= \min (r,\frac{1}{2})$, there exists a
$\tilde{s}>s$ such that \eqref{9-11-1} holds. Thus using standard
interpolation theory (see \cite{tr}) yields
$$\omega(x,t)\in L^\infty ([0, T^*); H^s(\Omega)).$$ The proof
of Lemma 3.1 is finished.

\vspace{3mm}

The following are some known facts, of which the proofs are omited
here.

\vspace{2mm}

\textbf{Lemma 3.2} (see \cite{bb})~~ For any $s\geq 1$ and $
1<p<\infty,$ there exists a positive constant $C$, depending only on
$\Omega, s, p,$ such that for any vector-valued function $w$,  one
has
$$\begin{array}{ll}\|w\|_{W^{s,p}(\Omega)} &\leq C[\|{\rm div}~w\|_{W^{s-1, p}(\Omega)} + \|{\rm
curl}~w\|_{W^{s-1, p}(\Omega)} \\[2mm] &\ \ \ \ + \|w\cdot \vec{n}\|_{W^{s - 1/p,
p}(\partial \Omega)} + \|w\|_{W^{s-1, p}(\Omega)}].\end{array}$$
Here $\vec{n}$ is the exterior normal vector on $\partial \Omega.$

\vspace{3mm} \textbf{Lemma 3.3}(Korn's Inequality)(see \cite{H})~~
Let $ \omega\in H^1(\Omega).$ Then there exists a constant $C$
depending only on the domain $\Omega$, such that
$$\|w\|_{H^1(\Omega)} \leq C(\|D(w)\|_{L^2(\Omega)} +
\|w\|_{L^2(\Omega)}).$$

Now we are ready to prove our main result.

\vspace{3mm} \textbf{Proof of Theorem 2.2}~~Let  $v^\nu = u^\nu -u.$
For any fixed $T < T^*$, one has for every $0<t\leq T$,
\begin{eqnarray}
&&\|u^\nu(t)\|_{L^2(\Omega)}^2 + 4\alpha \nu \int_0^t \int_{\partial
\Omega} |u^\nu|^2 + 4\nu \int_0^t \int_\Omega|D(u^\nu)|^2
\nonumber\\
&&\leq \|u^\nu(0)\|_{L^2(\Omega)}^2, \label{3.1}\\[3mm]
&& \|u(t)\|_{L^2(\Omega)}^2 = \|u_0\|_{L^2(\Omega)}^2. \label{3.1+}
 \end{eqnarray}
Here \eqref{3.1}  is the energy inequality for the leray weak
solution $u^\nu$ and \eqref{3.1+} is the energy equality for the
vortex patch solution $u$. Using $u$ as a test function in the weak
form satisfying by the Leray weak solution $u^\nu$ (see Definition
2.2), we obtain\begin{equation}\label{9-11-2}
\begin{array}{l}
\displaystyle\int_\Omega u^\nu \cdot u(t) dx
 + 2\alpha \nu \int_0^t \int_{\partial\Omega} u^\nu \cdot u
 dSd\tau
+ 2\nu \int_0^t \int_\Omega D(u^\nu) :D(u) dxd\tau\\[2mm]
+ \displaystyle\int_0^t \int_\Omega (u^\nu \cdot \nabla )u^\nu \cdot
u dxd\tau = \int_\Omega u^\nu(0)\cdot u_0 dx.
\end{array}\end{equation}
Adding \eqref{3.1} and \eqref{3.1+}  and then subtracting
\eqref{9-11-2},  one deduces
 \begin{eqnarray}\label{9-12-1}
 &&\ \ \ \ \frac{1}{2}
\|v^\nu(t)\|_{L^2(\Omega)}^2
+ 2\nu \int_0^t \int_\Omega |D(v^\nu)|^2 dxd\tau \nonumber\\[3mm]
&&\leq \frac{1}{2} \|v^\nu (0)\|_{L^2(\Omega)}^2 - \int_0^t
\int_\Omega (v^\nu \cdot \nabla)u \cdot v^\nu dxd\tau - 2\nu\alpha
\int_0^t
\int_{\partial \Omega} u^\nu \cdot v^\nu dSd\tau \nonumber\\[3mm]
&&\ \ \ \  - 2\nu \int_0^t \int_{\Omega} D(u): D(v^\nu) dxd\tau
\equiv \sum_{i=1}^{4} I_i.
\end{eqnarray}
 Now we estimate the terms on the right hand of \eqref{9-12-1}. By H\"{o}lder's
inequality,
\begin{eqnarray}\label{9-29-1}
|I_2|=\left|\int_0^t \int_\Omega (v^\nu \cdot \nabla)u \cdot v^\nu
dxd\tau\right| \leq \|\nabla u\|_{L^1(0, T; L^\infty(\Omega)}
\|v^\nu\|_{L^\infty(0,T; L^2(\Omega))}^2.
\end{eqnarray}
 For any
$0<\epsilon<\beta$, there exists $  0< s = \beta - \epsilon$ such
that $u\in L^\infty(0, T; H^s(\Omega))$. Using the duality between
$H^s(\Omega)$ and $H^{-s}(\Omega)$ (note that when $s< \frac{1}{2},
\ H^s(\Omega) = H_0^s(\Omega)$), one has
$$\begin{array}{ll}
& \left|\displaystyle\int_\Omega D(u): D(v^\nu)dxd\tau\right| \\[4mm] \leq &
C \|D(u)\|_{H^s(\Omega)} \|D(v^\nu)\|_{
H^{-s}(\Omega)}\\[2mm] \leq & C [\|{\rm curl}~u\|_{
H^s(\Omega)}+ \|u\|_{L^2(\Omega)}]\cdot \| v^\nu \|_{ H^{1-s}(\Omega)}\\[2mm]
\leq & C [\|\omega \|_{H^s(\Omega)} + \|u\|_{L^2(\Omega)}]\cdot
\|v^\nu\|_{L^2(\Omega)}^s
\|v^\nu\|_{H^1(\Omega)}^{1-s} \\
\leq & C [\|\omega \|_{H^s(\Omega)} + \|u\|_{L^2(\Omega)}] \|v^\nu
\|_{L^2 (\Omega)}^s (\|v^\nu \|_{L^2(\Omega)} + \|D(v^\nu)
\|_{L^2(\Omega)} )^{1-s},
\end{array}
$$
where the second inequality is the result of Lemma 3.2, the third
one is from an interpolation inequality, and the fourth one is due
to Lemma 3.3.

Since $v^\nu = u^\nu -u$, one has
$$
|I_3|=2\nu\left|\alpha \int_{\partial \Omega} u^\nu \cdot v^\nu
dS\right| \leq 2\nu\left|\alpha \int_{\partial \Omega} u \cdot v^\nu
dS\right| + 2\nu\left|\alpha \int_{\partial \Omega} |v^\nu|^2
dS\right|.
$$

Note that
$$\begin{array}{ll}
&\nu \left|\alpha \displaystyle\int_{\partial \Omega} u \cdot v^\nu dS\right| \\[3mm]
\leq & \nu \|\alpha \|_{L^\infty(\partial \Omega)}
\|u\|_{L^2(\partial
\Omega)} \|v^\nu\|_{L^2(\partial \Omega)}\\[2mm]
\leq & C\nu \|\alpha \|_{L^\infty(\partial \Omega)}
\|u\|_{H^{1/2}(\Omega)} \|v^\nu\|_{H^{1-s}(\Omega)} \\[2mm]
\leq & C\nu\|\alpha \|_{L^\infty(\partial
\Omega)}(\|\omega\|_{H^{s}(\Omega)} + \|u\|_{L^2(\Omega)})
\|v^\nu\|_{H^{1-s}(\Omega)}\\[2mm]
\leq & C\nu \|\alpha \|_{L^\infty(\partial \Omega)}
\|v^\nu\|_{L^2(\Omega)}^{\frac{2s}{1+s}}
+\frac{\nu}{4}\|D(v^\nu)\|_{L^2(\Omega)}^2,
\end{array}$$
and
$$\begin{array}{ll}
&\nu \left|\alpha \displaystyle\int_{\partial \Omega} |v^\nu|^2 dS\right|\\[3mm]
\leq &C \nu \|\alpha\|_{L^\infty(\partial \Omega)}
\|v^\nu\|_{H^{\frac12}(\Omega)}^2\\[2mm]
\leq & C\nu \|\alpha\|_{L^\infty(\partial
\Omega)}\|v^\nu\|_{L^2(\Omega)}(\|v^\nu\|_{L^2(\Omega)} +
\|D(v^\nu)\|_{L^2(\Omega)})\\[2mm]
\leq & C\nu \|\alpha\|_{L^\infty(\partial
\Omega)}\|v^\nu\|_{L^2(\Omega)}^2 +
\frac{\nu}{4}\|D(v^\nu)\|_{L^2(\Omega)}^2.
\end{array}
$$
The term $I_3$ is estimated as
\begin{eqnarray}\label{9-29-2}
I_3\le C\nu \|\alpha\|_{L^\infty(\partial
\Omega)}(\|v^\nu\|_{L^2(\Omega)}^2 +
\|v^\nu\|_{L^2(\Omega)}^{\frac{2s}{1+s}})+
\frac{\nu}{2}\|D(v^\nu)\|_{L^2(\Omega)}^2.
\end{eqnarray}

 Moreover, from  (\ref{3.1}), one deduces
$$\begin{array}{l}
\ \ \ \ \|u^\nu(t)\|_{L^2(\Omega)}^2 + 4\nu \int_0^t
\|D(u^\nu)\|_{L^2(\Omega)}^2 d\tau \\ [2mm]\leq
\|u^\nu(0)\|_{L^2(\Omega)}^2 +\nu \|\alpha\|_{L^\infty(\partial
\Omega)}\int_0^t \left( \|u^\nu\|_{L^2(\Omega)}
+\|D(u^\nu)\|_{L^2(\Omega)} \right)^{\frac{1}{2}}
\|u^\nu\|_{L^2(\Omega)}^{\frac12} d\tau\\[2mm]
\leq \|u^\nu(0)\|_{L^2(\Omega)}^2+  C\nu\int_0^t
\|u^\nu\|_{L^2(\Omega)}^2 d\tau+ 2\nu \int_0^t
\|D(u^\nu)\|_{L^2(\Omega)}^2 d\tau,
\end{array}$$ which implies that $\|u^\nu\|_{L^\infty(0,T;L^2(\Omega))}$ and
$\|v^\nu\|_{L^\infty(0, T; L^2(\Omega))}$ are uniformly bounded by
some constant $C$ depending on $M$, $T$ and
$\|\alpha\|_{L^\infty(\partial \Omega)}$. Hence by the Young's
inequality, since $\frac{2s}{s+1} < 1$,
\begin{eqnarray}\label{9-29-3}
&& |I_4|=\nu \left|\displaystyle\int_\Omega D(u): D(v^\nu)dxd\tau\right| \nonumber\\[4mm]
\leq && C\nu \|v^\nu\|_{L^2(\Omega)} + C\nu \|v^\nu
\|_{L^2(\Omega)}^s
\|D(v^\nu)\|_{L^2(\Omega)}^{1-s} \nonumber\\[2mm]
\leq && C\nu\|v^\nu\|_{L^2(\Omega)}^{\frac{2s}{1+s}} +\frac{\nu}{4}
\|D(v^\nu)\|_{L^2(\Omega)}^2,
\end{eqnarray}
where $C$ is a constant depending on $\Omega, s, M, T,
\|\alpha\|_{L^\infty(\partial \Omega)}$.

Putting \eqref{9-29-1}-\eqref{9-29-3} into \eqref{9-12-1},  we get
$$\begin{array}{l}\ \ \ \ \frac{1}{2}\|v^\nu(t)\|_{L^2(\Omega)}^2 \\\leq
\frac{1}{2}\|v^\nu(0)\|_{L^2(\Omega)}^2 + \displaystyle\int_0^t
(\|\nabla u\|_{L^\infty(\Omega)}+C)\|v^\nu\|_{L^2(\Omega)}^2 d\tau +
C\nu \int_0^t \|v^\nu \|_{L^2(\Omega)}^{\frac{2s}{1+s}} d\tau.
\end{array}$$ By the Gronwall lemma, we deduce that
$$\|v^\nu(t)\|_{L^2(\Omega)}^{\frac{2}{1+s}} \leq C\|v^\nu
(0)\|_{L^2(\Omega)}^{\frac{2}{1+s}} + C\nu t,$$ where $C$ is a
constant depending on $\epsilon, u, T,$
$\|\alpha\|_{L^\infty(\partial \Omega)}$ and $M\equiv \sup_\nu
\|u^\nu(0)\|_{L^2(\Omega)}$. The proof of the theorem is finished.

\vspace{2mm}

{\bf Acknowledgements:}   The authors would like to thank Professor
Zhouping Xin for his interest on this topic and  valuable
discussions. The paper was started when the first author was
visiting the Institute of Mathematical Sciences (IMS) of The Chinese
University of Hong Kong.


\vspace{2mm}
\begin{thebibliography}{100}
\bibitem{ad}Abidi, H., Danchin, R., Optimal bounds for the inviscid
limit of Navier-Sokes equations, {\it Asymptot. Anal.}, {\bf
38}(1)(2004), 35-46.

\bibitem{bc} Bertozzi, A.L., Constantin, P., Global regularity for
vortex patches, {\it Comm. Math. Phys.}, {\bf 152}(1) (1993), 19-28.

\bibitem{bb}Bourguignon J.P., Brezis H., Remarks on the Euler
Equation, {\it Jour. Func. Anal.}, {\bf 15}(1974), 341-363.

\bibitem{che1} Chemin, J.Y., Persistance de structures
g\'eom\'etriques dans les fluides incompressibles bidimensionnels,
{\it Ann. Sci. $\acute{E}$cole Norm. Sup. (4)}, {\bf 14}(2) (1993),
517-542.

\bibitem{che} Chemin, J.Y., A remark on the inviscid limit for
two-dimensional incompressible fluids, {\it Comm. Partial
Differential Equations}, {\bf 21}(11-12) (1996), 1771-1779.

\bibitem{cmr}Clopeau T., Mikeli\'{c} A. and Robert R., On the
vaninshing viscosity limit for the 2D incompressible Navier-Stokes
equations with the friction type boundary conditions, {\it
Nonlinearity}, {\bf 11} (1998), 1625-1636.

\bibitem{cw1} Constantin, P., Wu, J., Inviscid limit for vortex
patches, {\it Nonlinearity}, {\bf 8}(5) (1995), 735-742.

\bibitem{cw2} Constantin, P., Wu, J., The inviscid limit for
non-smooth vorticity, {\it Indiana Univ. Math. J.}, {\bf 45}(1)
(1996), 67-81.

\bibitem{de} Depauw, N., Poche de tourbillon pour Euler 2D dans un
ouvert $\grave{a}$ bord, {\it J. Math. Pures Appl.}, {\bf 78}(1995),
313-351.

\bibitem{du}Dutrifoy A., On 3-d vortex patches in bounded domains,
{\it Comm. P. D. E.}, {\bf 28}(2003), 1237-1263.

\bibitem{GS}Gamblin, P., Saint Raymond, X., On three-dimensional
vortex patches, {\it Bull. Soc. Math. France}, {\bf 123}(3)(1995),
375-424.

\bibitem{H}Horgan C. O., Korn's inequalities and their applications
in continuum mechanics, {\it SIAM Rev. }, {\bf 37}(4)(1995),
491-511.

\bibitem{ip}Iftimie, D., Planas, G., Inviscid limits for the
Navier-Stokes equations with Navier friction boundary conditions,
{\it Nonlinearity}, {\bf 19}(2006), 899-918.

\bibitem{JM} J$\ddot{a}$ger, W., Mikeli$\acute{c}$, A., On the
roughness-induced efective boundary conditions for an incompressible
viscous flow, {\it J. Differential Equations}, {\bf 170}(2001),
96-122.

\bibitem{JN} Jiu, Q.S., Niu, D.J., Vanishing viscous limits for the 2D lake equations
\\ with Navier boundary conditions, {\it J. Math. Anal. Appl.}, {\bf
338} (2008), 1070-1080.

\bibitem{ke}Kelliher, J. P., Navier-Stokes equations with Navier
boundary conditions for a bounded domain in the plane, {\it SIAM J.
Math. Anal.}, {\bf 38}(2006), 210-232.

\bibitem{P1}
Lions, P.L.,  Mathematical topics in fluid mechanics. {V}ol.
1:Incompressible models,
  Oxford Lecture Series in Mathematics and its Applications,  Oxford University Press, New York, 1996.

\bibitem{lnp}Lopes Filho, M. C., Nussenzveig Lopes, H. J., and Planas G., On the inviscid limit for 2D incompressible flow with
Navier friction condition, {\it SIAM J. Math. Anal.}, {\bf
36}(2005), 1130-1141.

\bibitem{mo}Masmoudi, N., Remarks about the inviscid limit of the Navier-Stokes
system, {\it Comm. Math. Phy.}, {\bf 270}(2007), 777-788.

\bibitem{na}Navier, C. L. M. H., Sur les lois de l'\'{e}quilibrie et du
mouvement des corps \'{e}lastiques, {\it Mem. Acad. R. Sci. Inst.
France}, {\bf 369}(1827).

\bibitem{S}Sueur, F., Vorticity internal transition layers for the
Navier-Stokes equations, {\it arXiv: 0812.2145v1}.

\bibitem{Swan} Swann, H. S. G., The convergence with vanishing
viscosity of nonstationary Navier-Stokes flow to ideal flow in
$R^3$, {\it Trans. Amer. Math. Soc.}, {\bf 157}(1971), 373-397.

\bibitem{tr}Triebel, H., Theory of function spaces, {\it Volume 100
of Monographs in Mathematics}, Birkh\"{a}user Verlag, Basel, (2006).

\bibitem{X1}
Xiao, Y. L., Xin, Z. P., On the vanishing viscosity limit for the 3D
Navier-Stokes equations  with a slip boundary condition, {\it Comm.
Pure Appl. Math.}, {\bf 60}(7) (2007), 1027-1055.


\end{thebibliography}
\end{document}